\documentclass[review]{elsarticle}
\usepackage{amsfonts,amssymb,amsmath,amsbsy,amsthm,amscd,bbm}
\usepackage[english]{babel}
\usepackage{graphicx}
\usepackage{color}
\usepackage{floatrow}
\newtheorem{theorem}{Theorem}[section]

\newtheorem{corollary}{Corollary}
\newtheorem{proposition}{Proposition}

\def\br{{\bf r}}
\renewcommand{\le}{\leqslant}
\renewcommand{\ge}{\geqslant}

\def\Div{\mbox{div}\,}

\def\bE{{\bf E}}

\def\dfrac#1#2{\displaystyle{#1\over #2}}

\def\bv{{\bf V}}
\def\bV{{\bf V}}

\begin{document}

\begin{frontmatter}

\title{Linearization method  and  sharp thresholds for
spherically symmetric multidimensional pressureless
Euler-Poisson equations
}

\author{Olga S. Rozanova}
\ead{rozanova@mech.math.msu.su} 
\author{Marko K. Turzynsky}
\ead{m13041@yandex.ru} 

\address[1]{ Moscow
State University, Leninskie Gory, Moscow 119991 Russia}
\address[2]{Russian University of Transport, Obraztsova, 9, Moscow 127055 and  Higher
School of Economics, Pokrovskiy Blvd, 11, Moscow, 109028, Russia}


\begin{abstract}
We show that the question about the criterion of a singularity formation for radially symmetric solutions to the Cauchy problem for a fairly wide class of equations related to the pressureless Euler-Poisson equations can be reduced to the study of solutions to a linear homogeneous ordinary differential equation.   In some cases, such a criterion can be obtained in terms of the initial data.  In the remaining cases, it is possible to construct a simple numerical procedure, on the basis of which the question about preserving smoothness for any set of initial data can be solved.
 \end{abstract}


\def\sign{\mathop{\rm sgn}\nolimits}





\end{frontmatter}

\section{Introduction}

The history of attempts to obtain criteria of a singularity formation for the solution of the Cauchy problem to the pressureless Euler-Poisson system is quite long. The pressureless Euler-Poisson system is interesting because it contains features of real physical models, but at the same time allows an accurate analytical study of threshold phenomena in terms of initial data, which occurs extremely rarely. The nature of the solution varies considerably depending on the assumptions made about the interaction force (attractive or repulsive) and the background density. For the case of one spatial variable, the question about the exact identification of the initial data corresponding to a globally smooth solution is almost completely solved in \cite{ELT}. However, the transfer of the results to the case of radially symmetric solutions turned out to be very difficult. A review can be found in \cite{Bhat23}, \cite{Tan}.

This work is inspired by the possibility of obtaining the criterion for a singularity formation in terms of initial data in the repulsive case in the exceptional spatial dimension 4 \cite{R_exept}.
The question arose as to what other classes of systems this method could be extended to. We show that the success in obtaining a criterion in terms of the initial data is related to the possibility of obtaining a first integral of some auxiliary system. In the multidimensional case, the criteria turn out to be quite cumbersome, since they involve not only derivatives of the initial data, but also the data themselves. Therefore, the space of initial data sets corresponding to a globally smooth solution is four-dimensional.

At the same time, such a criterion in terms of the existence of a zero of a particular solution of some linear homogeneous differential equation can be obtained for a wider class of systems. Although it does not give explicit conditions on the initial data under which the solution preserves global smoothness, it can be simply implemented numerically and allows the possibility of checking any initial data.

We are also substantially interested in the connection with the theory of linear differential equations, the results of which can be applied to the new area.

We study a class of extended Euler-Poisson equations
\begin{eqnarray}\label{EP}
\dfrac{\partial n }{\partial t} + \Div(n \bv)=0,\quad
\dfrac{\partial \bv }{\partial t} + \left( \bv \cdot \nabla \right)
\bv = -\,k  \nabla \Phi\, - \mu \bv +m \br , \quad -\Delta \Phi =n-c,
\end{eqnarray}
where the scalar functions $n$ (density),
$\Phi$ (force potential), and the vector $\bv$ (velocity) depend on
the time $t$ and the point $x\in {\mathbb R}^d $, $d\ge 1$, $ {\br}=(x_1, \dots ,x_d)$. Here
$c \ge 0$ is the density background, $\mu={\rm const}\ge 0$ is the friction coefficient, $m=\rm const$ is the intensity of the quadratic confinement, $k={\rm const}\ne 0$. The sign of $k$ corresponds to a repulsive (plus) or attractive (minus) force.

Many important models can be reduced to \eqref{EP}.
In particular,  it is the pressureless Euler-Poisson equations with a quadratic confinement \cite{Carrillo}
\begin{eqnarray}\label{CC}
\dfrac{\partial n }{\partial t} + \Div(n \bv)=0,\quad
\dfrac{\partial \bv }{\partial t} + \left( \bv \cdot \nabla \right)
\bv = -\int_{{\mathbb R}^d}\, {\nabla}_{ x} N(x-{ y}) n(t,{ y}) \, d{ y} -{\br},
\end{eqnarray}
where $-N({ x}) $ is the fundamental solution of the Laplace operator in $d$ - dimensional space, i.e. $-\Delta N ({x})=\delta({ x})$, $d\ge 2$. Indeed, the right hand side term in the second equation \eqref{CC} is $-\nabla \Psi$, where $\Psi $ is the solution of   $\Delta \Psi =n-d$. Thus, \eqref{CC} coincides with \eqref{EP} for $c=d$, $m=\mu=0$.

Further, the Euler-Poisson equations with a nonlocal pressure term (e.g. \cite{Brunelli})
\begin{eqnarray}\label{NL}
\dfrac{\partial n }{\partial t} + \Div(n \bv)=0,\quad
\dfrac{\partial \bv }{\partial t} + \left( \bv \cdot \nabla \right)
\bv = -\int_{{\mathbb R}^d}\, \nabla_{ x} N({ x}-{ y}) n(t,{ y}) \, d{ y},
\end{eqnarray}
similarly can be reduced to \eqref{EP} with $c=m=\mu=0$.

First, we write \eqref{EP} in a more convenient form
and introduce $\bE= \nabla \Phi$.
Under the assumption
 that the solution is sufficiently smooth and $\bE$ vanishes as $|\br|\to\infty$
  we
 obtain
\begin{eqnarray}\label{n}
n=c- \,\Div \bE,\label{n}\end{eqnarray} and remove $ n $
from \eqref{EP}. The result is
\begin{eqnarray}\label{4}
\dfrac{\partial \bv }{\partial t} + \left( \bv \cdot \nabla \right)
\bv = \, - k \bE -\mu \bv + m \br,\quad \frac{\partial \bE }{\partial t} + \bv \Div \bE
 = {c} \bV.
\end{eqnarray}

Denote $ r=|{\bf r}|$, and consider radially symmetric solutions   depending only on
$r$,
\begin{eqnarray}\label{sol_form}
\bv=F(t,r) \,\br,\quad \bE=G(t,r)\, \br,\quad
n=n(t,r)\quad n_0=n_0(r).
\end{eqnarray}

 Consider the initial data
\begin{equation}\label{CD1}
(\bv, \bE) |_{t=0}= (F_0(r) \, {\bf r}, G_0(r) \, {\bf r} ), \quad (F_0(r)
, G_0(r) ) \in C^2(\bar {\mathbb R}_+),
\end{equation}
such that  $n|_{t=0}\ge 0$.

For the local in $t$ well-posedness of the  Cauchy problem \eqref{EP}, \eqref{CD1} we refer to \cite{Tan}. Notice that the formation of singularity is associated with infinite gradient of the solution.

We call a solution of \eqref{4},  \eqref{CD1} smooth for $t\in
[0,t_*)$, $t_*\le\infty$, if the functions $F$ and $G$ in
\eqref{sol_form} belong to the class $ C^1([0,t_*)\times
\bar{\mathbb R}_+)$. The blow-up of solution implies that the
derivatives of solution tends to infinity as $t\to t_*<\infty$.

The paper is organized as follows.Sections \ref{Sec2} and \ref{Sec3} are devoted to deriving the equations of behavior of the solution components and their derivatives along the characteristic. In Section \ref{Sec4}, we perform a procedure of linearizing the equations for the derivatives and prove Theorem \ref{T1} stating that the question about the criterion for a singularity formation  in a fairly general situation can be reduced to studying the possibility of some solution of a linear equation vanishing. Section \ref{Sec5} classifies singular points on the phase plane associated with the solution components  depending on the assumptions about  the acting force ($k>0$ or $k<0$) and the background density ($c>0$ or $c=0$). In Section \ref{Sec6}, we consider the repulsive case $k>0$, $c>0$, the results are a  generalization of \cite{R_exept}. Further results concern dimension $d\ge 3$. In Section \ref{Sec7}, we present a general scheme for reducing the question of the possibility of obtaining the desired criterion to the possibility of obtaining an analytical solution of some linear homogeneous second-order equation. In Section \ref{Sec8}, we apply this scheme to the case $k>0$, $c=0$ and obtain a criterion for the formation of a singularity for the case of zero initial velocity (Theorem \ref{Tn00}). As a consequence, we obtain a similar criterion  for a system of gas dynamics equations with nonlocal pressure. The case  $d=4$ is again different from the others. Namely, we show that it is possible to obtain a criterion for arbitrary initial data, which is most simple for zero initial velocity. In Section \ref{Sec9} we analyze the case $k<0$, $c>0$ and show that analytical results can be obtained only for $d=4$ (in terms of special functions), Theorem \ref{Tk-1}. If analytical results cannot be obtained, numerical illustrations are given to confirm that the qualitative structure of the boundary of the sets of initial data that guarantee global smoothness is the same as in the case where analytical results are possible.
  Moreover, we provide a sufficient condition for the singularity formation for all $d\ge 1$.
Section \ref{Sec10} is devoted to a discussion about the possibility of applying the described technique to other physically significant situations.

\bigskip

\section{The behavior along characteristics}\label{Sec2}

A great advantage of the problem in the radially symmetric case is that the behavior of all components of the solution and their derivatives is completely described by their behavior along the Lagrangian trajectory, i.e. along the only characteristic of the first-order quasilinear system. This allows the use of a well-developed technique for studying systems of ordinary differential equations.

First of all, we note that from \eqref{4}, \eqref{CD1} it follows that the functions $F$ and $G$ satisfy the following Cauchy problem:
\begin{equation}\label{fgsys}
\frac{\partial G}{\partial t}+Fr\frac{\partial G}{\partial r}={c} F - {d}FG,\quad
\frac{\partial F}{\partial t}+Fr\frac{\partial F}{\partial r}=-F^2- kG - {\mu} F - m,
\end{equation}
\begin{equation*}\label{fgcon}
(F(0,r),\,G(0,r))=(F_0(r),\,G_0(r)),\quad (F_0(r),\,G_0(r))\in C^2(\overline{\mathbb{R}}_+)
\end{equation*}

Along the characteristic
\begin{equation}\label{ch}
\dot r=Fr,
\end{equation}
which starts from the point $r_0\in\left[0,+\infty\right)$, the system \eqref{fgsys}  takes the form
\begin{equation}\label{sysch}
\dot G ={c} F - {d}FG,\quad \dot F = -F^2-m-kG-{\mu} F.
\end{equation}

Let us introduce the notation
\begin{eqnarray*}
&&\mathcal{D}=\Div{\bf V}, \quad  \lambda=\Div{\bf E}, \\ &&J_{ij}=\partial_{x_i} V_i \partial_{x_j} V_j
-\partial_{x_j} V_i \partial_{x_j} V_i,\quad  i\neq j, \quad J=\sum\limits_{i,j=1, i<j}^{d} J_{ij}.
 \end{eqnarray*}
 The number of terms in the sum for $J$
is equal to $\frac{{d}({d}-1)}{2}$. It is easy to show that the relations
 \begin{equation*}\label{conn}
 \mathcal{D}={d}F+F_r r,\quad \lambda={d} G+G_r r,\quad J=({d}-1)F F_r r+\frac{{d}({d}-1)}{2} F^2,
 \end{equation*}
are satisfied, which implies
\begin{equation}\label{j}
2J=2({d}-1)\mathcal{D}F-({d}-1){d}F^2.
\end{equation}
Then from  system \eqref{4} we have
\begin{eqnarray*}\label{sys1}
&&\frac{\partial \mathcal{D}}{\partial t}+({\bf V} \cdot \nabla\mathcal{D})=-\mathcal{D}^2+2J-k\lambda-m{d}-\mu\mathcal{D},\\
&&\frac{\partial \lambda}{\partial t}+({\bf V} \cdot \nabla\lambda)=\mathcal{D}({c}-\lambda)+(\nabla c,{\bf V}).
\end{eqnarray*}
Thus, taking into account equality \eqref{j} along the characteristic outgoing from the point $r_0$, we obtain
\begin{eqnarray}\label{sys2}
&&\mathcal{\dot D}=-\mathcal{D}^2+2({d}-1)F\mathcal{D}-k\lambda-m{d}-({d}-1){d}F^2-\mu\mathcal{D},\\\nonumber &&\dot\lambda=\mathcal{D}({c}-\lambda)+c' r F.
\end{eqnarray}
This is a quadratically nonlinear system with $F$, $r$ and $c(r)$ as  coefficients, which can be found from \eqref{ch}, \eqref{sysch}. We can also consider   \eqref{sysch} and \eqref{sys2} as a quadratically nonlinear system of 4 equations on $G$, $F$, $\mathcal{D}$, $\lambda$, in which system \eqref{sysch} is closed. Introducing new variables $u=\mathcal{D}-{d}F$, $v=\lambda -{d}G$, we obtain
\begin{eqnarray}\label{sys3}
&&\dot u = -u^2-(2F+\mu)u-kv-{d}(m+\mu F),\\\nonumber &&\dot v = -uv + ({c}-{d}G)u-{d}Fv+{c'} r F.
\end{eqnarray}

\section{Linearization of  system \eqref{sys3}}\label{Sec3}

The main tool for further research is the following theorem (\cite{Radon}, \cite{Riccati}).

\begin{theorem} [The Radon lemma (1927)]
\label{T2} A matrix Riccati equation
\begin{equation}
\label{Ric}
 \dot W =M_{21}(t) +M_{22}(t)  W - W M_{11}(t) - W M_{12}(t) W,
\end{equation}
 {\rm (}$W=W(t)$ is a matrix $(n\times m)$, $M_{21}$ is a matrix $(n\times m)$, $M_{22}$ is a matrix  $(m\times m)$, $M_{11}$ is a matrix  $(n\times n)$, $M_{12} $ is a matrix $(m\times n)${\rm )} is equivalent to the homogeneous linear matrix equation
\begin{equation}
\label{Lin}
 \dot {\mathcal Y} =M(t) {\mathcal Y}, \quad M=\left(\begin{array}{cc}M_{11}
 & M_{12}\\ M_{21}
 & M_{22}
  \end{array}\right),
\end{equation}
 {\rm (}${\mathcal Y}={\mathcal Y}(t)$  is a matrix $(n\times (n+m))$, $M$ is a matrix $((n+m)\times (n+m))$ {\rm )} in the following sense.

Let on some interval ${\mathcal J} \in \mathbb R$ the
matrix-function $\,{\mathcal Y}(t)=\left(\begin{array}{c}\mathfrak{Q}(t)\\ \mathfrak{P}(t)
  \end{array}\right)$ {\rm (}$\mathfrak{Q}$  is a matrix $(n\times n)$, $\mathfrak{P}$  is a matrix $(n\times m)${\rm ) } be a solution of \eqref{Lin}
  with the initial data
  \begin{equation*}\label{LinID}
  {\mathcal Y}(0)=\left(\begin{array}{c}I\\ W_0
  \end{array}\right)
  \end{equation*}
   {\rm (}$ I $ is the identity matrix $(n\times n)$, $W_0$ is a constant matrix $(n\times m)${\rm ) } and  $\det \mathfrak{Q}\ne 0$ on ${\mathcal J}$.
  Then
{\bf $ W(t)=\mathfrak{P}(t) \mathfrak{Q}^{-1}(t)$} is the solution of \eqref{Ric} with
$W(0)=W_0$ on ${\mathcal J}$.
\end{theorem}

To apply this theorem, we rewrite system \eqref{sys3} in matrix form \eqref{Ric}. Then $W=(u\,\,v)^T$,
\begin{eqnarray*}
  M_{11}=\left(\begin{matrix}  0   \end{matrix}\right), &&  M_{12}=\left(\begin{matrix} 1 & 0   \end{matrix}\right), \\
 M_{21}= \left(\begin{matrix} -{d}(m+\mu F)  \\ {c'} r F  \end{matrix}\right), &&  M_{22}=
\left(\begin{matrix} -2F-\mu & -k \\{c} - {d}G & -{d} F \end{matrix}\right).
\end{eqnarray*}

\noindent Thus, we obtain the Cauchy problem for a system linear with respect to the functions $q, p_1, p_2$:
\begin{eqnarray}\label{2}&&
\left(\begin{array}{ccc} \dot q \\ \dot p_1 \\ \dot p_2\end{array}\right)=\left(\begin{array}{ccc} 0 & 1 & 0\\
-{d}(m+\mu F)& -2F-\mu & -k\\c' r F & {c}-{d}G &-{d}F \end{array}\right) \left(\begin{array}{ccr} q \\ p_1 \\ p_2\end{array}\right),\\ \label{2CD}&& \left(\begin{array}{ccr} q(0) \\ p_1(0) \\ p_2(0)\end{array}\right)=\left(\begin{array}{ccr} 1 \\ u_0 \\ v_0\end{array}\right),
\end{eqnarray}
in which the coefficients $r(t), G(t), F(t)$ are found from
\eqref{ch}, \eqref{sysch}.
System \eqref{2} in all cases can be reduced to a single third-order equation
\begin{eqnarray}\nonumber
\dddot q &+&[(2+{d})F+\mu]\ddot q+ [2({d}-1)F(F+\mu)-kG({d}+2)+
 m({d}-2)+k{c(r)}]\dot q\\&+&[\mu {d} ({d}-1)F^2 +(md^2-\mu^2 d +kc'r)F-\mu{d}(m+kG)]q=0.\label{qdif}
\end{eqnarray}

\section{Criterion of the singularity formation in terms of auxiliary function}\label{Sec4}
We see from Theorem \ref{T2} that the derivatives of radially symmetric solutions of  system \eqref{4} go to infinity when the auxiliary function $q(t)$ vanishes at some point of the semiaxis $t>0$. Thus, we can formulate the following theorem.

\begin{theorem}\label{T1} Suppose that the components of the radially symmetric solution $(\bE, \bv)$ of  system \eqref{4} have at most linear growth as $r\to\infty$. Then the solution of the Cauchy problem \eqref{4}, \eqref{CD1} preserves smoothness for all $t>0$ if and only if $q(t)$, the component of the solution of the Cauchy problem for the linear system \eqref{2}, does not vanish for all $t>0$.
\end{theorem}

Indeed, if we restrict ourselves to solutions that grow in space no faster than a linear function, we consider bounded $(F,G)$, i.e. bounded coefficients of the system \eqref{2}. Note that solutions that grow in space as a linear function play an important role in the theory of constructing solutions to the Euler-Poisson equations \cite{R22_Rad}, \cite{RT}. In particular, explicit solutions can be constructed in this class.

Theorem \ref{T1} is implicit, and it is difficult to construct on its basis a set of initial data for problem \eqref{4}, \eqref{CD1} corresponding to a smooth solution. However, it gives a very simple numerical algorithm that allows one to check whether any given initial data \eqref{CD1} belong to this class. Examples of the application of this algorithm to particular cases of  system \eqref{4} can be found in \cite{R22_Rad}, \cite{Roz_doping}, \cite{RD24}, we will not dwell on the numerical results now, but will try to identify cases when Theorem \ref{T1} allows us to obtain as a consequence a criterion for preserving smoothness in terms of the initial data.

Success is ensured by the possibility of obtaining the first integrals of  system \eqref{ch}, \eqref{sysch}, \eqref{2}.

\section{Phase curves of the system \eqref{sysch} and classification of equilibria}\label{Sec5}
In what follows, due to the desire to obtain an analytical first integral of the system \eqref{sysch},
we  restrict ourselves to a constant value of $c$, as well as a zero value of $\mu$.  If an analytical first integral cannot be obtained, then it is still possible to obtain various estimates of the solution that lead to some sufficient conditions on the initial data that ensure the global smoothness of the solution (e.g. \cite{RD24}), but a criterion cannot be obtained in this way.

Thus,
for $\mu=0$
on the phase plane $(G,\,F)$ the phase curves of system \eqref{sysch} are symmetric with respect to the axis $F=0$, and the system itself is reduced to one equation
\begin{equation*}\label{one}
\frac{1}{2}\frac{dF^2}{dG}=\frac{F^2+kG+m}{{d} G-{c}},
\end{equation*}
which is linear with respect to $F^2$. For ${d}=2$ its direct integration yields
 \begin{eqnarray}\label{fsim}
 &2F^2=(2G-c)k\ln|2G-c|-kc-2m+C_2 (2G-c),\\\nonumber
 &C_2=\frac{2F_0^2(r_0)+kc+2m}{2G_0(r_0)-c}-k\ln|2G_0(r_0)-c|.
 \end{eqnarray}

 For ${d}=1$ or ${d}\ge 3$ we similarly have
  \begin{eqnarray}\label{fsim1}
 &F^2=C_{ d} |{ d} G-c|^{\frac{2}{\bf d}}+\frac{k(2G-c)}{{ d}-2}-m,\\
 \label{Cd}&C_{ d}=\frac{(F_0^2(r_0)+m)({ d}-2)-k(2G_0(r_0)-c)}{({ d}-2)|{ d} G_0(r_0) - c|^{\frac{2}{ d}}}.
 \end{eqnarray}

\begin{proposition}  On the phase plane $(F,\,G)$ corresponding to  system \eqref{sysch}, the motion occurs either in the half-plane $G<\frac{c}{ d}$, or in the half-plane $G>\frac{c}{ d}$, or on the line $G=\frac{c}{ d}$.
\end{proposition}
\proof
 From \eqref{sysch} taking into account   \eqref{ch} we obtain
\begin{equation*}
\frac{{d}G}{F(c-{d}  G)}=\frac{dr}{Fr}
\end{equation*}
and
\begin{equation*}
c - {d}  G=\mbox{const} \cdot r^{-{d}}.
\end{equation*}
From here we see that the sign of the expression $(c-{d}  G)$ does not change, i.e.
\begin{equation*}
\mbox{sign}(c-{d}  G)=\mbox{sign}(c-{d}  G_0(r_0)).
\end{equation*}
$\Box$


The system \eqref{sysch} for $\mu=0$ has the following equilibria:
\begin{enumerate}
\item For $m{d}+ck> 0$ there is a unique equilibrium  ($F=0$, $G=-\frac{m}{k}$),  a {\bf center};

 \item For  $m{d}+ck=0$ there is a unique equilibrium  ($F=0$, $G=-\frac{m}{k}$),   a {\bf saddle-node};
\item For $m{d}+ck<0$ there are three points:

 ($F=0$, $G=-\frac{m}{k}$), which is a {\bf saddle},

 ($F=\sqrt{-m-\frac{ k c}{ d}}$, $G=\frac{c}{ d}$) --- a {\bf stable node}, and ($F=-\sqrt{-m-\frac{k c}{d}}$,

 $G=\frac{c}{ d}$) --- an {\bf unstable node}.

\end{enumerate}
Note that the change
\begin{equation}\label{change}
G_1=G+\frac{m}{k}, \quad c_{1}= c+ \frac{d m}{k}
\end{equation}
reduces the equilibrium  ($F=0$, $G=-\frac{m}{k}$) to zero, but at the same time, if no additional conditions are imposed, it may turn out that the new background value $c_{1}$ is negative. In what follows we restrict ourselves to the case $c_{1}\ge 0$, since otherwise we change $G_1$ to $-G_1$ and
$k$ to $-k$.

\begin{figure}[htb]
\begin{minipage}{0.3\columnwidth}
\includegraphics[scale=0.2]{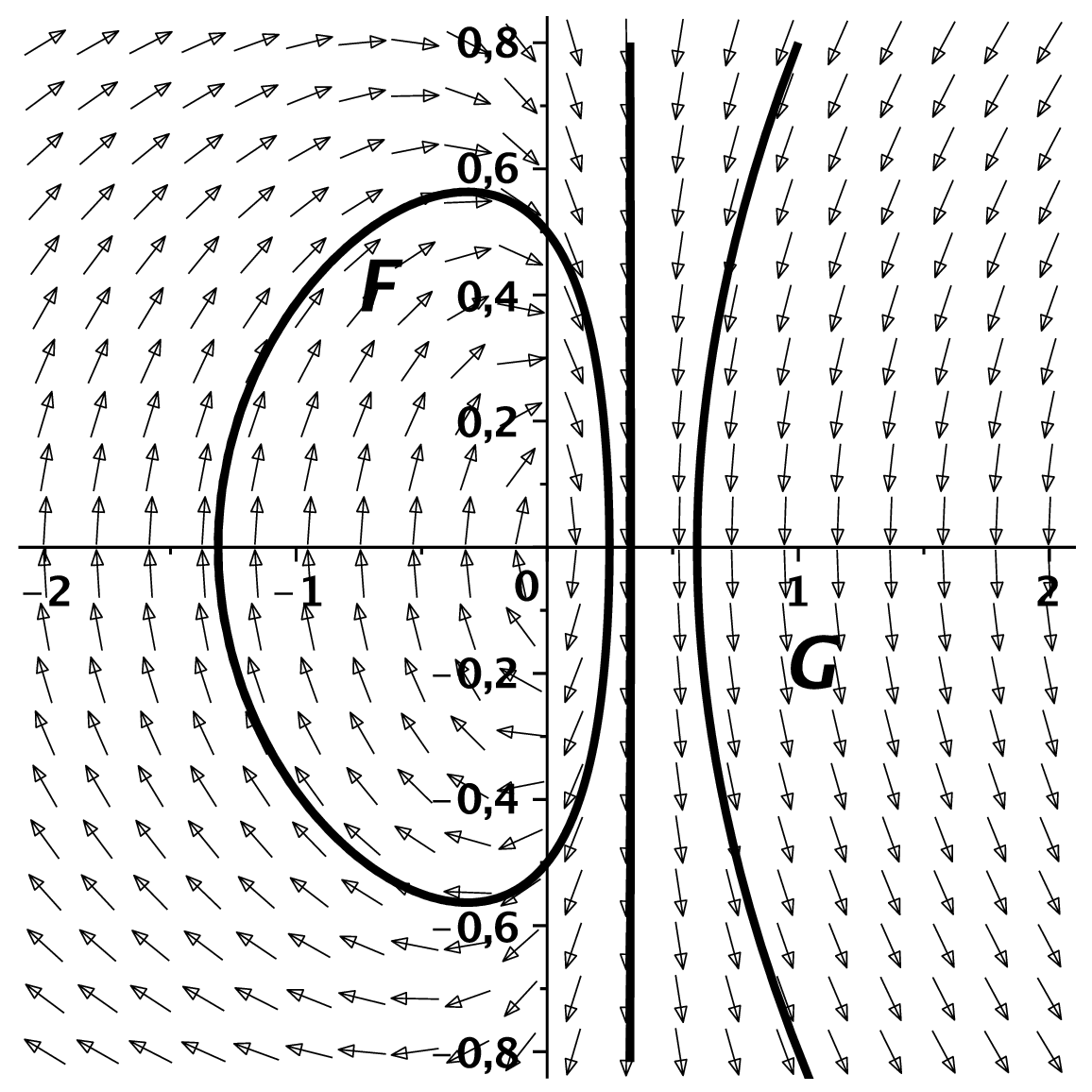}
\end{minipage}
\hspace{0.3cm}
\begin{minipage}{0.3\columnwidth}
\includegraphics[scale=0.2]{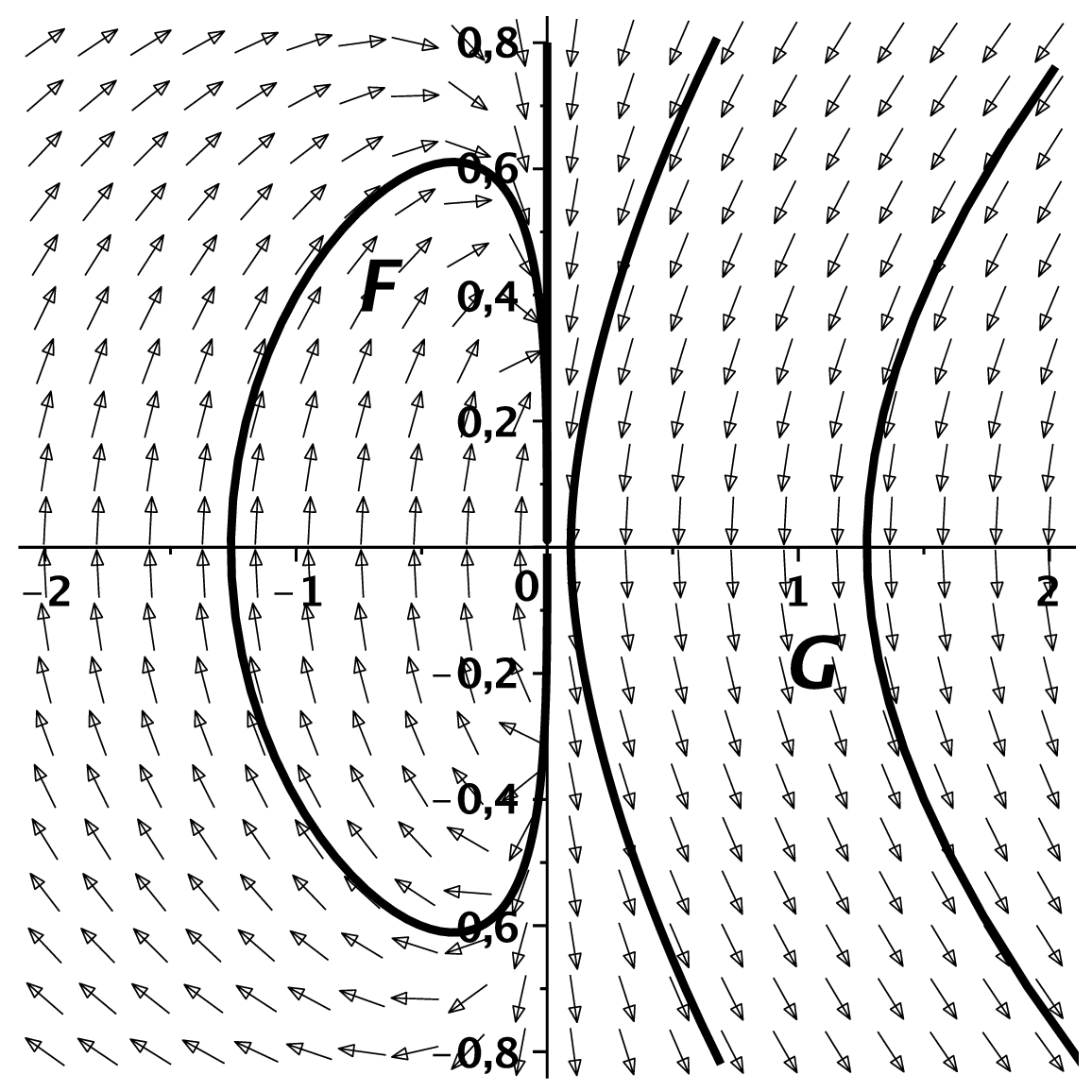}
\end{minipage}
\hspace{0.3cm}
\begin{minipage}{0.3\columnwidth}
\includegraphics[scale=0.2]{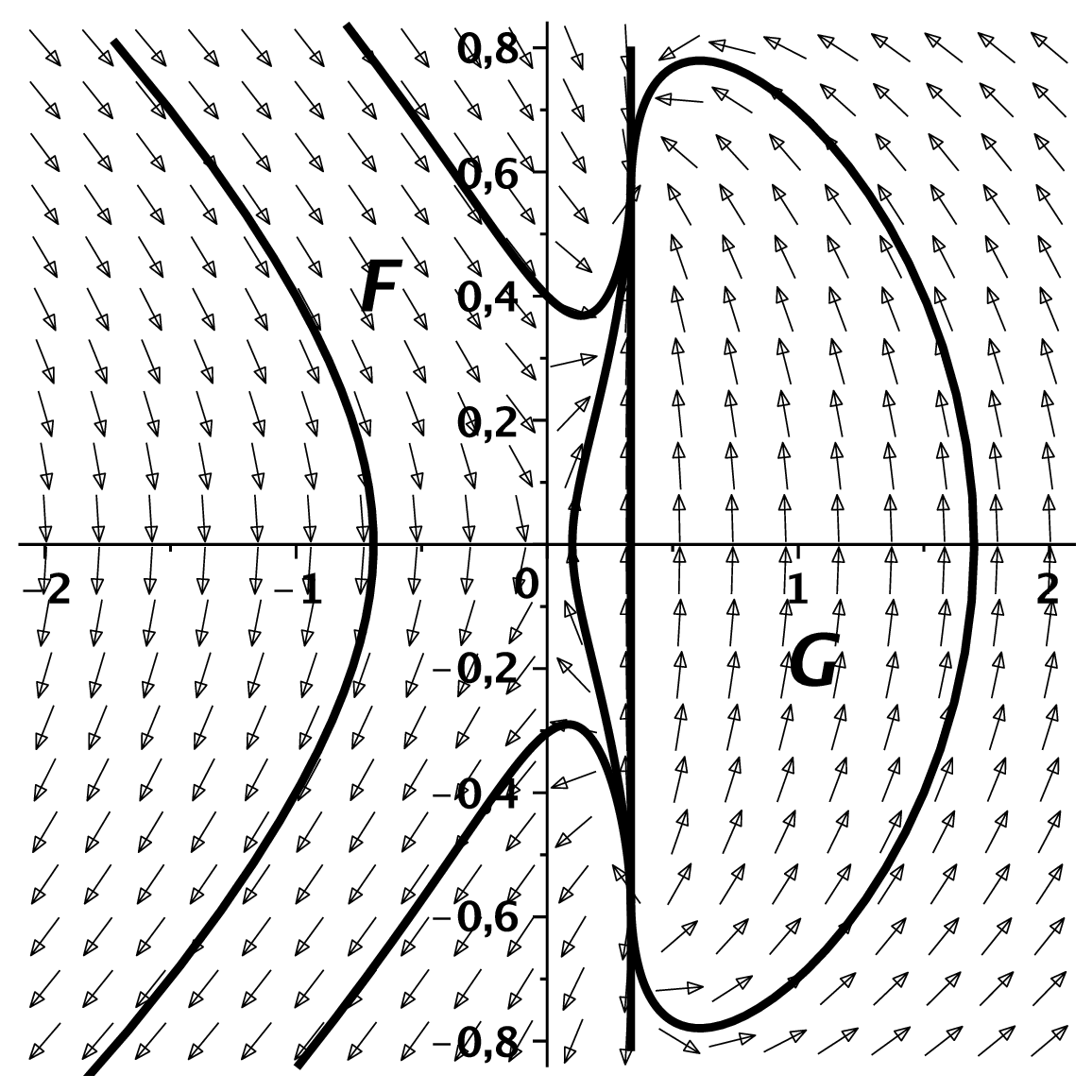}
\end{minipage}

\caption{The behavior of solution of \eqref{sysch} for $d=3$. Left: oscillatory case, $k=1$, $c=1$, $m=0$. Center: oscillatory case I, $k=1$, $c=0$, $m=0$. Right: non-oscillatory case II, $k=-1$, $c=1$, $m=0$.}\label{Pic1}
\end{figure}

We will always consider only the motion in the half-plane $G<\frac{c}{ d}$.  Indeed, since the positivity of  density and \eqref{n} imply
 \begin{equation}\label{npos}
    {\rm div} \bE = rG_r +d G < c
 \end{equation}
  and we assumed that $G_r$ is bounded in zero, then for $r_0=0$ we have $G<\frac{c}{ d}$.  For $r>0$ (here $r$ is the independent variable), from \eqref{npos} we have
\begin{equation}\label{n+}
r G(t,r) < -(d-1) \int\limits_0^r G(t, \xi) \xi d\xi +cr.
\end{equation}
 Assume that for a point $G_0(r_0)>\frac{c}{ d}$, then  $G>\frac{c}{ d}$ for all $t$ along this characteristic. For $k>0$
the function $G(t)$ tends to plus infinity as $t\to t_*<0$ for all possible data, therefore \eqref{n+} leads to a contradiction with \eqref{npos}.
For $k<0$, $c=0$,  \eqref{n+} cannot be valid for $G>0$. For $k<0$, $c>0$ the domain $G>\frac{c}{ d}$ does not contain the origin, therefore all solutions corresponding to this domain are unbounded (for example, the affine solutions $(F(t) \br, G(t) \br)$, for which $G_r=0$ and $G>\frac{c}{ d}$ implies $n<0$).

 \section{Oscillatory case (the equilibrium is a center), $k>0$, $c>0$.}\label{Sec6}
For $m{d}+ k c> 0$ the equilibrium  ($F=0$, $G=-\frac{m}{k}$) is a center. Let us perform the change of variables \eqref{change}
and therefore consider the equilibrium   ($F=0$, $G_1=0$) and $c_{1}>0$. To avoid cluttering the notation, we  omit the index 1.

 Thus,  there is a neighborhood of the origin $(0,0)$  in which the phase trajectories on the plane $(F, G)$ are closed and the solution with initial data lying on these trajectories is periodic. Let us study the question when all phase trajectories lying in the half-plane $G<\frac{c}{ d}$ are closed.

We will prove the following proposition.
\begin{proposition}\label{P1} Phase trajectory of system \eqref{sysch} starting from $r_0\in \overline{\mathbb R}_+$ for $k>0$ are closed in the half-plane $G<\frac{c}{ d}$

1. for $d=1 $
if and only if $F_0^2< {c} - 2 G_0$;

2. for $d\ge 2 $ for any initial point from  $G<\frac{c}{ d}$.

In the half-plane $G>\frac{c}{ d}$ all trajectories are unbounded.
\end{proposition}

\proof

1. For $d=1$, $k>0$, system \eqref{sys2} coincides with \eqref{sysch}, and the inequality $F_0(r_0)^2< {c} - 2 G_0(r_0)$ implies the known criterion for the global smoothness of a solution to the Cauchy problem $\bv'_0(r_0)^2< {c} - 2 \bE_0'(r_0)$, $r_0\in \mathbb R$, see \cite{RChZAMP21}.

2.  For $d\ge 2$, $k>0$, the boundedness of all trajectory for $G<\frac{c}{d}$ is proved in \cite{R22_Rad}, Lemma 2. The unboundedness for $G<\frac{c}{ d}$ follows similarly from the comparison of degrees of $G$ in the right hand sides of \eqref{fsim} and \eqref{fsim1}.
 $\Box$

The respective phase trajectories are presented in Pic.1, left.

For the case of periodic solutions $(F,G)$, the following result applies \cite{Roz_doping}, \cite{Carrillo}: if the period of oscillations $(F,G)$ depends on the initial point of the trajectory (the equilibrium is not an isochronous center), then the Lagrangian trajectories corresponding to different initial points of the characteristic necessarily intersect within a finite time, and, accordingly, the solution of the Cauchy problem \eqref{4}, \eqref{CD1} loses smoothness within a finite time. Therefore, the only globally smooth solution corresponds to the equilibrium position itself, that is,
$\bv=0$, $\bE=-\frac{m}{k}\br$.

Note that the isochronicity of the system means that it has an additional first integral. In the theory of Hamiltonian systems, such a situation is called superintegrability. Moreover, since Pouncar\'e, it has been known that there is a transformation by which an isochronous system can be linearized \cite{Romanovski}.

For the study of isochronous oscillations in our case, the following Sabatini criterion is convenient \cite{Sabatini}.
\begin{theorem}\label{S}
Let us consider a  Li\'enard type equation
\begin{equation}\label{Lienard}
\ddot y+ f (y) \dot y+g(y)=0,
\end{equation}
where $f, g$  are analytic, $g$ odd, $f (0)=g(0)=0$, $g'(0)>0$. Then
$\mathcal O =(y,\dot y)=(0,0)$ is a center if and only if $f$ is odd and
$\mathcal O )$ is an isochronous center if and only if
\begin{equation}\label{tau}
\tau(y) :=\left(\int\limits_0^y sf(s) ds \right)^2-y^3 (g(y)-g'(0)y)=0.
\end{equation}
\end{theorem}

Note that  system \eqref{sysch} can be rewritten as \eqref{Lienard}:
\begin{equation*}
\ddot F+(2+{d})F\dot F+(c k+m{d})F+{d}F^3=0.
\end{equation*}
Having calculated $\tau(F)$, as \eqref{tau}, we see that
in the case $(c k+m{d})>0$ the system has an isochronous center if and only if $(2+{d})^2=9{d}$, i.e. ${d}=1$ and ${d}=4$.

Thus, only in dimensions ${d}=1$ and ${d}=4$ there exists an open neighborhood of the equilibrium position in the $C^1$ norm, the initial data from which correspond to a globally smooth solution.

For the equilibrium $(0,0)$ and $c>0$ there are results concerning the possibility of constructing globally smooth solutions. Namely, for $d\ne 1$, $d\ne 4$ the only possibility of this kind is a simple wave $F=F(G)$ \cite{R22_Rad}, \cite{R_exept}, for $d= 1$ the conditions that single out the initial data corresponding to smooth solutions are found in \cite{ELT}, \cite{RChZAMP21}, for $d= 4$ the criterion for the formation of singularities in terms of the initial data is found in \cite{R_exept}. By applying the substitution \eqref{change} we can shift the equilibrium and apply the already known results.

\section{  Behavior of derivatives, formal scheme}\label{Sec7}

Let us study the behavior of derivatives for the case $d\ge 3$. To this aim we introduce a new variable
$M=|c-d G|^\frac{2}{d}$, $M>0$, therefore
\begin{equation*}\label{GF}
 G=\frac{c-M^\frac{d}{2}}{d}\qquad F^2(M)=C_d M-\frac{2 k}{d (d-2)} \, M^\frac{d}{2}- \frac{k c}{d}-m,
\end{equation*}
$C_d$ is defined in \eqref{Cd}. Let us denote $P(t)=p_1(t)$ and $R(t)=\dot P(t)$.
Thus, from the first equation of \eqref{sysch}, \eqref{fsim1}, and \eqref{qdif} we get
\begin{eqnarray}\label{qs}
  \dfrac{dq}{dM} &=& -\frac{P}{2 F M} \\
   \dfrac{dP}{dM} &=&  -\frac{R}{2 F M}\label{p1s} \\
   \dfrac{dR}{dM} &=& \frac{1}{2 F M} (Q_2(M) R +Q_1(M) P + Q_0(M) q),\label{p2s}
\end{eqnarray}
where
\begin{eqnarray*}
& & Q_0(M) = d^2 m F(M), \\ &&Q_1(M)= 2(d-1) F^2(M)-k (d+2) G(M) + m(d-2) + k c,\\&& Q_2(M)=(d+2) F(M).
\end{eqnarray*}

If we assume that the shift  \eqref{change} is made, then $Q_0(M)=0$,
and $P(M)$ is a solution of the second order linear homogeneous ODE
\begin{eqnarray}\label{Y}
&&  Y''(M)+ S_1(M) Y'(M) + S_2(M)  Y(M)=0, \\
&&  S_1(M)=\frac{-2 d (d-1)(d-2)C_d M + kd ((d-2)d c - 2 M^\frac{d}{2})}{4d (d-2) C_d M^2 - 2k (c d(d-2) + 4 M^\frac{d}{2}) M},\nonumber\\
&&   S_2(M)=\frac{2 d (d-1)(d-2)C_d M - c kd (d-2)  M^\frac{d}{2} -c k(d-2)(d^2-d+2)}{4d (d-2) C_d M^3 - 2k (c d(d-2) + 4 M^\frac{d}{2}) M^2}.\nonumber
\end{eqnarray}
If equation \eqref{Y} is solved subject to initial conditions
 \begin{eqnarray*}\label{YIC}
  &&Y(M_0)=P(M_0)=u_0, \\  &&Y'(M_0)=-2 F(M_0) M_0 R(M_0)=2 F(M_0) M_0 (2 F(M_0) u_0+k v_0),
 \end{eqnarray*}
 (see \eqref{2}) where ($M_0=|c-d G_0|^\frac{2}{d}$), then
\begin{equation}\label{qM}
  q(M)=1-\int\limits_{M_0}^M \frac{Y(\xi)}{2 \xi F(\xi)}\, d\xi,
\end{equation}
and then we have to analyze this function. For a globally smooth solution it should be positive for all possible $M$.

Of course, this program can hardly be performed in the general case. In what follows we show that sometime this problem still can be solved.

\section{ Non-oscillatory case I (the equilibrium is a saddle-node),  $c=0$.}\label{Sec8}

Again we assume that the shift  \eqref{change} is made, therefore $c=0$. Further, we notice that we can restrict ourselves to the case $k>0$, since if $k<0$ (the attractive case) in the half-plane $(G, F)$ all trajectories are unbounded and both components tend to minus-infinity as $t\to t_*<\infty$.

Thus, the only equilibrium is the origin $(0,0)$. Exactly as in the case of
Proposition \ref{P1}, we prove that for $d\ge 2$ all trajectories are bounded in the half-plane $G<0$.
For $d=1$ the solution $(G,F)$ is bounded in the half-plane $G<0$ if and only if $F\ge 0$ or $F_0(r_0)^2<  - 2 G_0(r_0)$, $F<0$ for every $r_0$ and for every initial data (see also \cite{T24}).

   In the half-plane $G>0$ every solution is unbounded. See Pic.1, center.

    The attractive case can be reduced to the repulsive case if  we change $k$ to $-k$ and $G$ to $-G$. Thus, in the attractive case in the physical half-plane $G>0$   every non-trivial solution blows up.  It is quite natural since in the pressureless case in does not exist a force that counteracts the gravity. The presence of pressure changes the situation dramatically (e.g. \cite{Brenner}).

   The plan of this section is the following. First we obtain the criterion of the singularity formation for the initial conditions \eqref{CD1} with $\bv_0(r)=0$, which implies $F_0=v_0=0$, since the formulation of the result is too cumbersome. Then for the case
   $d=4$, where the solution can be expressed in a compact algebraic form, we obtain the criterion in the general case.

   Below we set $k=1$.

   \bigskip

   \subsection{$d\ge 3$, $\bv_0=0$}\label{8.1}

   Equation \eqref{Y} has following solutions for $F_0$:
   \begin{eqnarray*}
     Y_1(M)=M \sqrt{M_0^\frac{d-2}{2}- M^\frac{d-2}{2}}\,& {\mathcal F}\left(\left[-\frac{1}{2},\frac{d-1}{d-2}\right], \left[\frac{1}{d-2}\right], \left(\frac{M}{M_0}\right)^{\frac{d-2}{2}}\right), \\
     Y_2(M)=M^{\frac{d-1}{2}} \sqrt{M_0^\frac{d-2}{2}- M^\frac{d-2}{2}}\, & {\mathcal F}\left(\left[\frac{d-4}{2(d-2)},2\right], \left[2-\frac{1}{d-2}\right], \left(\frac{M}{M_0}\right)^{\frac{d-2}{2}}\right).
   \end{eqnarray*}
   Here ${\mathcal F}\left(\left[a,b\right], \left[c\right], z\right) $ is the Gaussian hypergeometric function \cite{hyper}, the solution of the equation
   \begin{equation*}
     z(1-z) \dfrac{d^2 w}{d z^2}+(c-(a+b+1)z) \dfrac{dw}{dz}-ab w=0,
   \end{equation*}
   $0<M\le M_0$, from the point $(M_0, 0)$ on the phase plane $(M,F)$ the trajectory moves from $M_0$ to zero. Further, $ Y_1(0)= Y_2(0)=0$, the limits of $ Y_1(M)$ and $Y_2(M)$ as $M\to M_0-0$ are finite.
   Moreover, one of the functions can be simplified as
   \begin{equation*}
    Y_1= \left(\frac{M}{M_0}\right)  \left(d \left(\frac{M}{M_0}\right)^\frac{d-2}{2}-2\right).
   \end{equation*}
   For $d\ge 4$ both solutions $ Y_1(M)$ and $ Y_2(M)$ are linearly independent and can be taken as a fundamental system, however,  $Y_1=Y_2$ for $d=3$, and we can add to the fundamental system together with $Y_1$ another solution
    \begin{equation*} \label{Y2_3}
 Y_2(M)=M \sqrt{\sqrt{M_0}- \sqrt{M}}\quad  {\mathcal F}\left(\left[-\frac{1}{2},2\right], \left[\frac{3}{2}\right], 1-\sqrt{\frac{M}{M_0}}\right).
   \end{equation*}

   Notice that for $d=3$ and $d=4$ the function $Y_2(M)$ is such that $Y_2(M_0)=0$. For $d\ge 5$ we introduce another fundamental system,
   $Y_1(M)$ and $\bar Y_2(M)=Y_2(M)-\frac{Y_2(M_0)}{Y_1(M_0)} Y_1(M)$, such that $\bar Y_2(M_0)=0$. For $d=3$ and $d=4$ we denote $\bar Y_2=Y_2$. Further,
    \begin{eqnarray*}
   &&p_1(M) = C_1 Y_1+C_2 \bar Y_2 , \\ &&p_2(M)= 2 M F(M) (C_1 Y'_1+C_2 \bar Y'_2)- 2 F(M)  (C_1 Y_1+C_2 \bar Y_2),
   \end{eqnarray*}
   $ u_0= p_1(M_0)=0= C_1 Y_1(M_0)$, therefore $C_1=0$, and
   \begin{equation}\label{C2}
    C_2 = \lim_{M\to M_0}\frac{v_0}{2 M_0 F(M) \bar Y_2'(M)}.
   \end{equation}
   Thus, according to \eqref{qM}, we obtain the following result.
   \begin{theorem}\label{Tn00}
   The solution to the problem \eqref{4}, \eqref{CD1} for $k=1$, $c=0$, $\bv_0=0$ preserves the initial smoothness if and only if
   for every $r\in \overline{\mathbb R}_+$ the following inequality holds:
   \begin{equation}\label{qcrit}
  1+C_2 \int\limits_{0}^{M_0} \frac{\bar Y_2(\xi)}{2 \xi F(\xi)}\, d\xi>0,
\end{equation}
where $C_2$ is defined in \eqref{C2}.
   \end{theorem}

   Notice that  $\bar Y_2'(M_0)=-\infty$, but since $F>0$ when $M_0$ tends to zero and $F(M_0)=F_0=0$, then it can be shown that the limit \eqref{C2} exists and the sign of $C_2$ coincides with the sign of $v_0$. Therefore if \eqref{qcrit} is not satisfied, then $v_0$ is negative.

   As a corollary of this theorem we can obtain the criterion for a singularity formation for the Euler-Poisson equations with the nonlocal pressure term \eqref{CC} with zero initial velocity.
   \begin{corollary} Let $d\ge 3$. The solution to the Cauchy problem \eqref{NL},
   \begin{equation*}
    (n, \bv)|_{t=0}=(n_0(r)>0, \bv=0),
   \end{equation*}
   keeps the initial smoothness for all $t>0$ if and only if condition \eqref{qcrit} holds with
   \begin{eqnarray*}
 M_0(r) = d|G_0(r)|^\frac{2}{d}, \quad G_0(r)=\frac{1}{r}\,\int_{{\mathbb R}^d}\, \nabla_{ x} N({ x}-{ y}) n_0({ y}) \, d{ y},
 \quad v_0(r)= n_0(r)-d G_0(r) 
     \end{eqnarray*}
     for every $r\in \mathbb R_+$.
     \end{corollary}
     \proof It is enough to notice that \eqref{CC} can be reduced to \eqref{EP} for $c=m=\mu=0$, $k=1$, for $\bE = \nabla  \Phi=\int_{{\mathbb R}^d}\, \nabla_{ x} N({ x}-{ y}) n({t, y}) \, d{ y}$. Further, ${\rm div} \bE =n = rG_r+ dG =v+dG$. $\Box$

   \subsection{$d=4$}

1.  First of all we obtain a corollary of Theorem \ref{Tn0} for $d=4$.

\begin{corollary}\label{Tn0}
   If $d=4$, then  the solution to the problem \eqref{4}, \eqref{CD1} for $k=1$, $c=0$, $\bv_0=0$ preserves the initial smoothness if and only if
   for every $r\in \overline{\mathbb R}_+$ the following inequality holds:
   \begin{equation}\label{qcrit4}
  1+\frac{2 v_0}{M_0^2} >0.
\end{equation}
 \end{corollary}
 \proof
  If $F_0=0$, then  $Y_1=M (2 M-M_0)$ and $Y_2=M^\frac{3}{2}\sqrt{M-M_0}$. It can be easily shown that $\frac{\bar Y_2(M)}{2 M F(M)}=1$, $C_2=\frac{2 v_0}{M_0^3}$ and
  the function $q(M)$ can be explicitly found as
  \begin{equation*}
    q(M)=-\frac{2 v_0}{M_0^3} (M-M_0) 
    +1,
  \end{equation*}
  and \eqref{qcrit} reads as \eqref{qcrit4}. $\Box$

  \bigskip

 2. For $\bv\ne 0$ the expression for $q(M)$ is more complicated but can still explicitly found by solving
\eqref{qs} - \eqref{p2s}:
\begin{eqnarray*}
 && q(M) = 1 + A_1(M_0, F_0; M) u_0 + A_2(M_0, F_0; M) v_0, \\
  &&M \in (0, M_+), \quad M_+=\frac{4F_0^2+M_0^2}{M_0}.
\end{eqnarray*}
The dependence of $A_1$ and $A_2$ on $M$ is quite cumbersome (it is the reason why we do not write out the coefficients here), however, it can be easily obtained by means of a computer algebra package.
Moreover,
\begin{equation*}
  q''(M)= \frac{B(M_0, F_0, u_0, v_0)}{(M (4 F_0^2+M_0^2-M_0 M))^\frac{3}{2}},
\end{equation*}
where $B$ depends only on the initial data. Therefore, the convexity of the graph of $q(M)$ preserves for all $ M \in (0, M_+)$. To check whether specific initial data $M_0, F_0, u_0, v_0$ belong to the smoothness domain in 4D space, it is necessary to investigate the positivity of $q(M)$ using standard methods (on $(0, M_0)$ for $F_0\ge 0$ or on $(0, M_+)$ for $F_0< 0$).

\section{ Non-oscillatory case II (the stable equilibrium is a node), $k<0$, $c>0$.}\label{Sec9}

First of all we notice that in contract to the previous cases the zero equilibrium is nor stable here. The asymptotically stable equilibrium on the plane $(F,G)$ in the point $(F_*=\sqrt{\frac{-k c}{d}}, G_*={\frac{ c}{d}} )$
(see Pic.1, right). This point corresponds to the affine solution
\begin{equation*}
  (\bV, \bE) = (F_* \br, G_* \br ),
\end{equation*}
which us unbounded as $r\to\infty$.

\subsection{Sufficient condition for the blow-up}

  We see that if the initial data are such that $(F_0,G_0)$ lie below than a separatrix of the saddle in the origin, then $F$ and $G$ tend to $-\infty$ and the solution blows up in a finite time (the fact that this time is finite follows from \eqref{fsim1} and \eqref{sysch}). 

The equation for the separatrices can be found from \eqref{fsim1}, choosing the constant $C_d$ such that this phase curve goes through the origin (see Fig.2): \begin{eqnarray}\label{separ}
  F^2 =   \frac{k c}{d-2}\,\left(\left(1-\frac{d G}{c}\right)^\frac{2}{d}+ \frac{2G}{c}-1\right).
   \end{eqnarray}

\begin{figure}[htb!]
\floatbox[{\capbeside\thisfloatsetup{capbesideposition={right,top},capbesidewidth=4cm}}]{figure}[\FBwidth]
{\caption{Separatrices of saddle equilibrium and the direction field for system \eqref{sysch} at $d=3$, $c=0$, $k=-1$, $G<\frac{c}{d}$. The blow-up region is below the separatrix emerging from the unstable node.}\label{Pic2}}
{\includegraphics[scale=0.2]{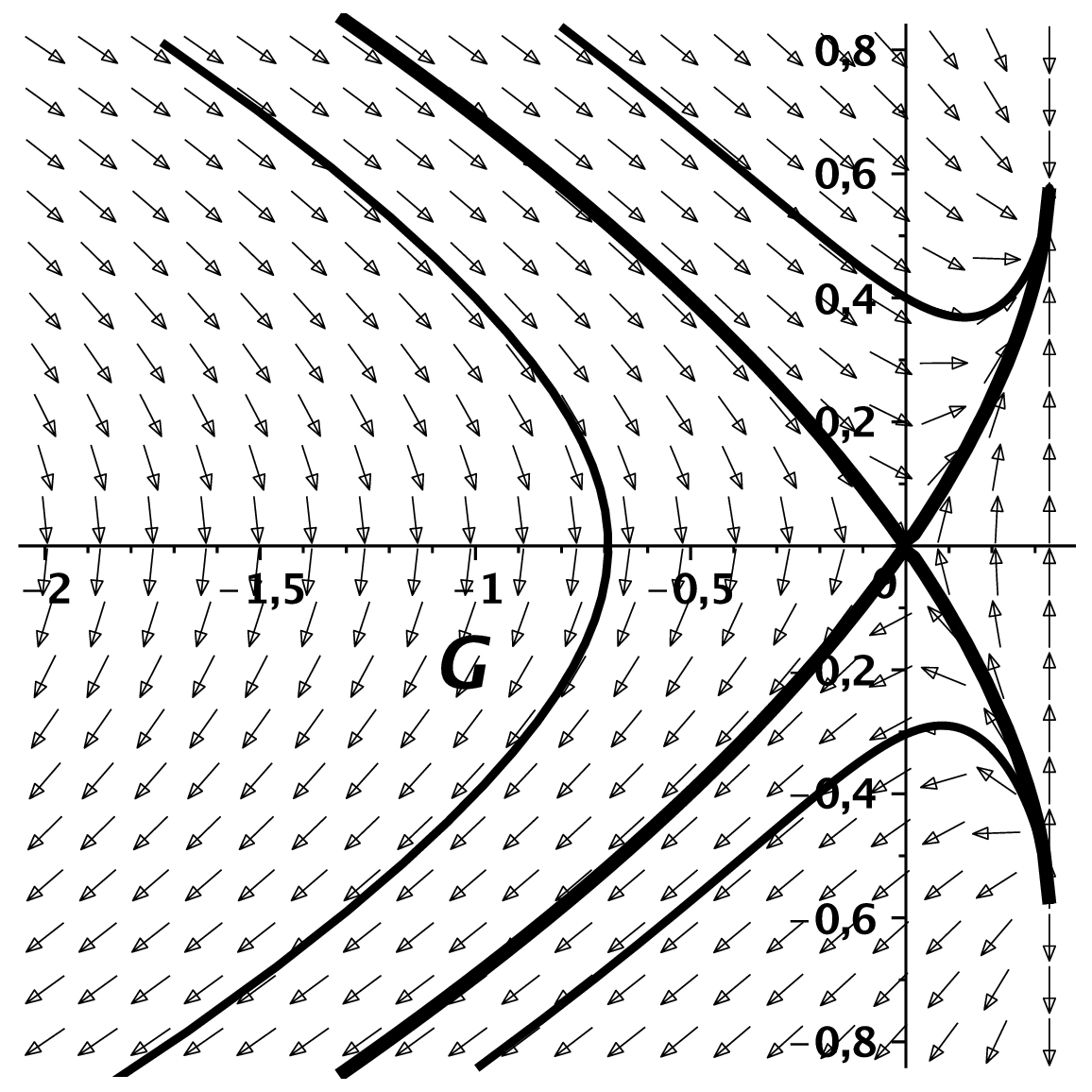}}
\end{figure}

   Thus, we can formulate the following statement.
   \begin{proposition}\label{Psep}
   If the initial data \eqref{CD1} are such that there exists $r_0\in \overline {\mathbb {R}}_+$ for which
   \begin{eqnarray}\label{separ_cond}
  F_0 <   - {\rm sign} G_0\, \sqrt{\frac{k c}{d-2}\,\left(\left(1-\frac{d G_0}{c}\right)^\frac{2}{d}+ \frac{2G_0}{c}-1\right)},
   \end{eqnarray}
   then the solution to the Cauchy problem \eqref{4}, \eqref{CD1} for $k<0$ and $c>0$ blows up in a finite time.
    \end{proposition}

  Notice that for $d=1$ condition \eqref{separ_cond} gives the criterion of the singularity formation, whereas for $d>1$ for the case
   of initial data are such that the inequality  opposite to  \eqref{separ_cond} is satisfied, we have to analyse the behavior of $u$ and $v$.
For $d=2$ the equation for separatrices, analogous to \eqref{separ} and the sufficient condition of the singularity formation can also be obtained from \eqref{fsim}.

\subsection{$d=4$}

The success of  the method of reducing to
\eqref{Y}  is ensured by the possibility of expressing its solution in known functions. It is not possible for an arbitrary $d$.
Only for $d=4$ the fundamental system of \eqref{Y}  can be expressed in the derivatives of the Heun functions \cite{Heun}.

As in Sec.\ref{8.1} we consider initial data with $\bv_0$, therefore $F_0=u_0=0$. For the sake of simplicity we set $c=1$, $k=-1$. The variable $M\in (0,1)$.

The fundamental system in this case consists of
\begin{eqnarray*}
&&Y_1(M)=\sqrt{2-M M_0} \, M^{\frac{3+\sqrt{2}}{2}}\,{\mathcal H}\left([a_1,q_1],[\alpha_1,\beta_1,\gamma_1,\delta_1];\frac{M}{M_0} \right), \\
&&Y_2(M)=\sqrt{2-M M_0} \, M^{\frac{3-\sqrt{2}}{2}}\,{\mathcal H}\left([a_2,q_2],[\alpha_2,\beta_2,\gamma_2,\delta_2]; \frac{M}{M_0}\right),
\end{eqnarray*}
where $,{\mathcal H}$ is Heun function, the solution of the problem
\begin{equation*}
     \dfrac{d^2 w}{d z^2}+\left(\frac{\gamma}{z}+\frac{\delta}{z-1}+\frac{\epsilon}{z-a}\right) \dfrac{dw}{dz}
     +\frac{\alpha \beta z -q}{z(z-1)(z-2)} w=0,\quad w(0)=1,\, w'(0)=\frac{q}{\gamma a},
   \end{equation*}
   and
\begin{eqnarray*}
&&  a_1= \frac{2}{M_0^2}, \, q_1= 1+ \frac{3 \sqrt{2}}{4} +\frac{2+\sqrt{2}}{2 M_0^2},\, \alpha_1=\frac{\sqrt{2}}{2}, \,\beta_1= 2+\frac{\sqrt{2}}{2},\, \gamma_1= 1+\sqrt{2}, \, \delta_1=\frac{1}{2} \\
  &&a_2= \frac{2}{M_0^2}, \, q_2= 1- \frac{3 \sqrt{2}}{4} +\frac{2-\sqrt{2}}{2 M_0^2},\, \alpha_2=-\frac{\sqrt{2}}{2}, \,\beta_2= 2-\frac{\sqrt{2}}{2},\, \gamma_2= 1-\sqrt{2}, \, \delta_2=\frac{1}{2}.
\end{eqnarray*}

One can check that $Y_1(0)=Y_2(0)=0$, $Y_1(M_0)\ne 0$, $Y_2(M_0)\ne 0$. As in Sec.\ref{8.1} we introduce another fundamental system,
   $Y_1(M)$ and
   \begin{equation}\label{Y-k}
   \bar Y_2(M)=Y_2(M)-\frac{Y_2(M_0)}{Y_1(M_0)} Y_1(M),
   \end{equation}
    such that $\bar Y_2(M_0)=0$.  Then
    \begin{eqnarray*}
   &&p_1(M) = C_1 Y_1+C_2 \bar Y_2 , \\ &&p_2(M)= 2 M F(M) (C_1 Y'_1+C_2 \bar Y'_2)- 2 F(M)  (C_1 Y_1+C_2 \bar Y_2),
   \end{eqnarray*}
   $ u_0= p_1(M_0)=0= C_1 Y_1(M_0)$, therefore $C_1=0$, and $C_2$ can be computed as \eqref{C2}.
    Thus, according to \eqref{qM}, we obtain the following analog of Theorem \ref{Tn00}.
   \begin{theorem}\label{Tk-1}
   The solution to the problem \eqref{4}, \eqref{CD1} for $d=4$, $k=-1$, $c=1$, $\bv_0=0$ preserves the initial smoothness if and only if
   for every $r\in \overline{\mathbb R}_+$ the following inequality holds:
   \begin{equation}\label{qqcrit}
  1+C_2 \int\limits_{0}^{M_0} \frac{\bar Y_2(\xi)}{2 \xi F(\xi)}\, d\xi>0,
\end{equation}
where $\bar Y_2$  and $C_2$ are defined in \eqref{Y-k} and \eqref{C2}, respectively.
   \end{theorem}

\subsection{$d\ge 3$, numerical result}
Since we are interested in the structure of the set of smoothness for any $d$, we can perform a numerical procedure based directly on
\eqref{2},  \eqref{sysch} with data \eqref{2CD} and $(F_0, G_0)$.

Some idea of how the set of initial  data is structured for a solution to be globally smooth can be obtained by considering the case of the equilibrium itself, when $F=F_*=\sqrt{-\frac{kc}{d}}  $,  $G=G_*=\frac{c}{d}  $. Then system \eqref{2} has constant coefficients and the solution can be explicitly found. Thus, if in a point $r_0\ge 0$ the initial data for  system \eqref{2}, \eqref{sysch}  are such that $(F_0=F_*, G_0=G_*)$, then $q(t)$ preserves positivity if and only if $u_0$ and $v_0$ satisfy
\begin{equation*}
1-\frac{d F_* u_0+ v_0}{1- d (G_*+2 F_*^2)}>0.
\end{equation*}

One can expect that in the general case the dependence in the criterion of the singularities formation on $u_0$ and $v_0$ is also linear. Indeed, let us consider  the system on $p_1$ and $p_2$ (a part of \eqref{2}). The coefficients $F(t)$ and $G(t)$ with the data above the separatrix \eqref{separ_cond} tends asymptotically to  $F_*$ and $G_*$ as $t\to \infty$.
 As follows from Theorem 8, Ch.2 \cite{Bellman}, the fundamental system consists of two solutions, having as $t\to \infty$ the same asymptotics that this system with constant coefficients $F_*$ and $G_*$, i.e. ${\rm const}\, e^{-\sqrt{cd}\, t}$ and ${\rm const}\, e^{-2\sqrt{\frac{c}{d}}\, t}$. Nevertheless, to find $q(t)$, we have to integrate $p_1$, and the initial $F_0$ and $G_0$ play the role.

 The numerical computations show that if we fix any two parameters in the quadruple $(F_0, G_0, u_0, v_0)$, the relationship between the rest two parameters on the borderline is (very close to) linear. Fig.3 presents this kind of pictures for $d=3$, $c=1$.

\begin{figure}[htb]
\begin{minipage}{0.3\columnwidth}
\includegraphics[scale=0.2]{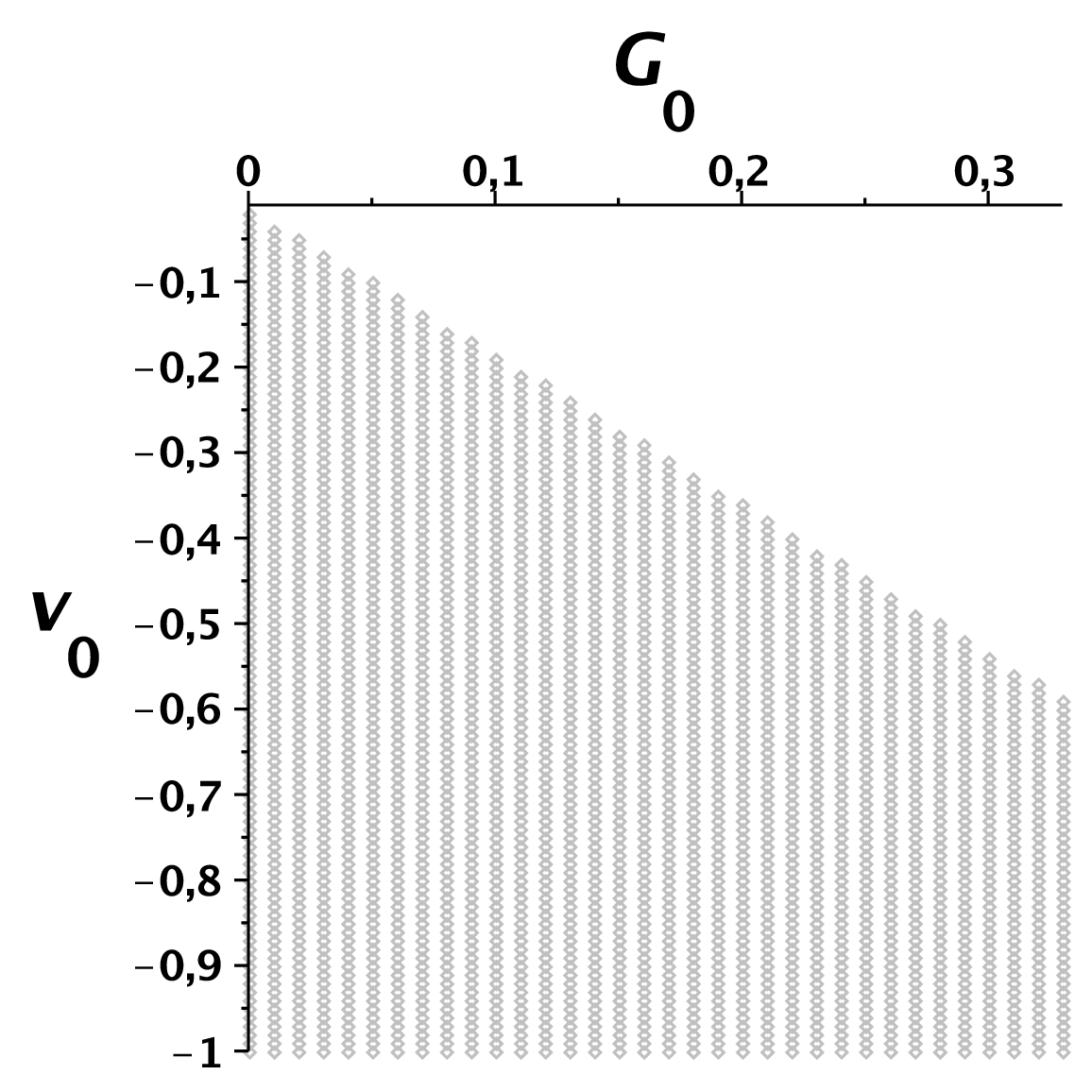}
\end{minipage}
\hspace{0.3cm}
\begin{minipage}{0.3\columnwidth}
\includegraphics[scale=0.2]{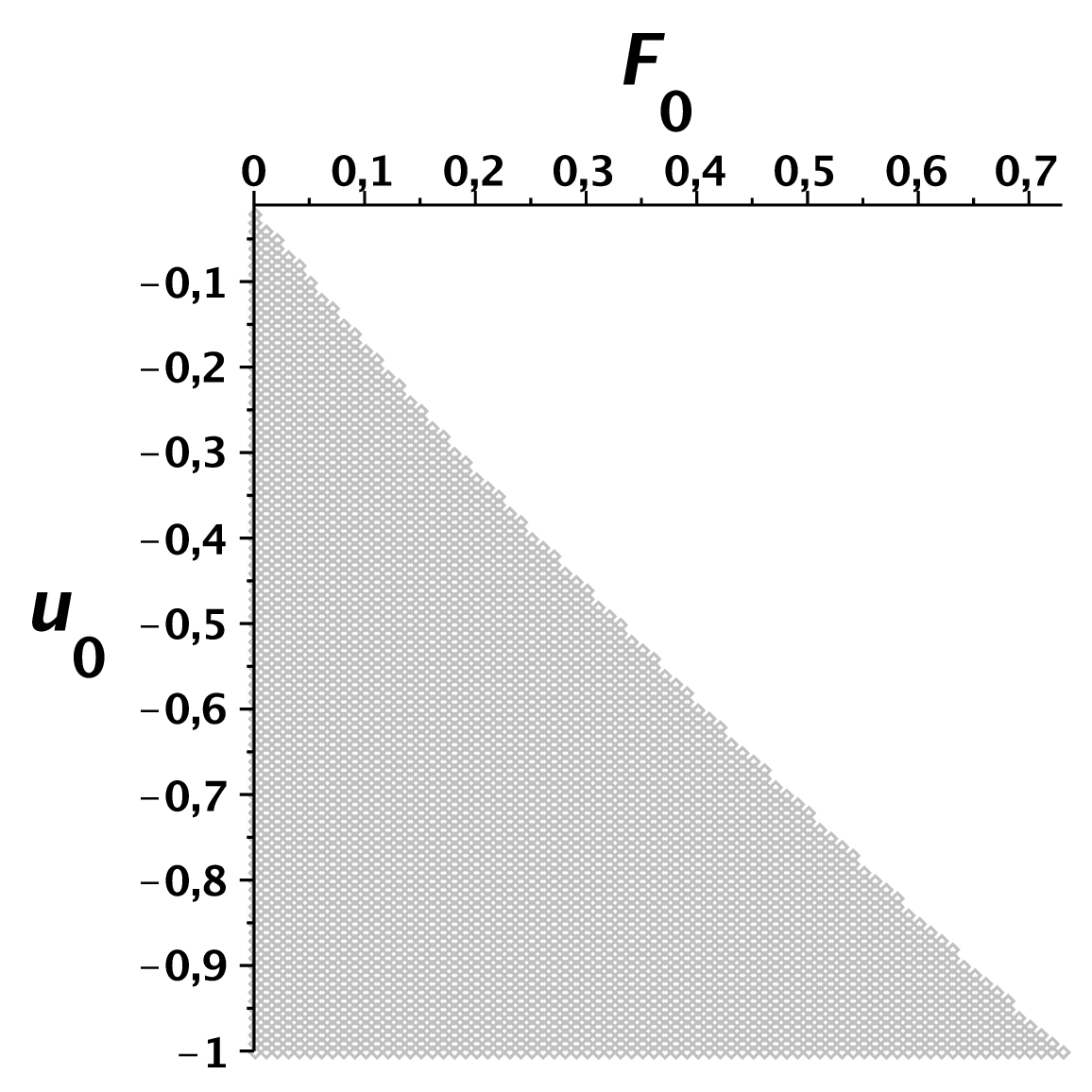}
\end{minipage}
\hspace{0.3cm}
\begin{minipage}{0.3\columnwidth}
\includegraphics[scale=0.2]{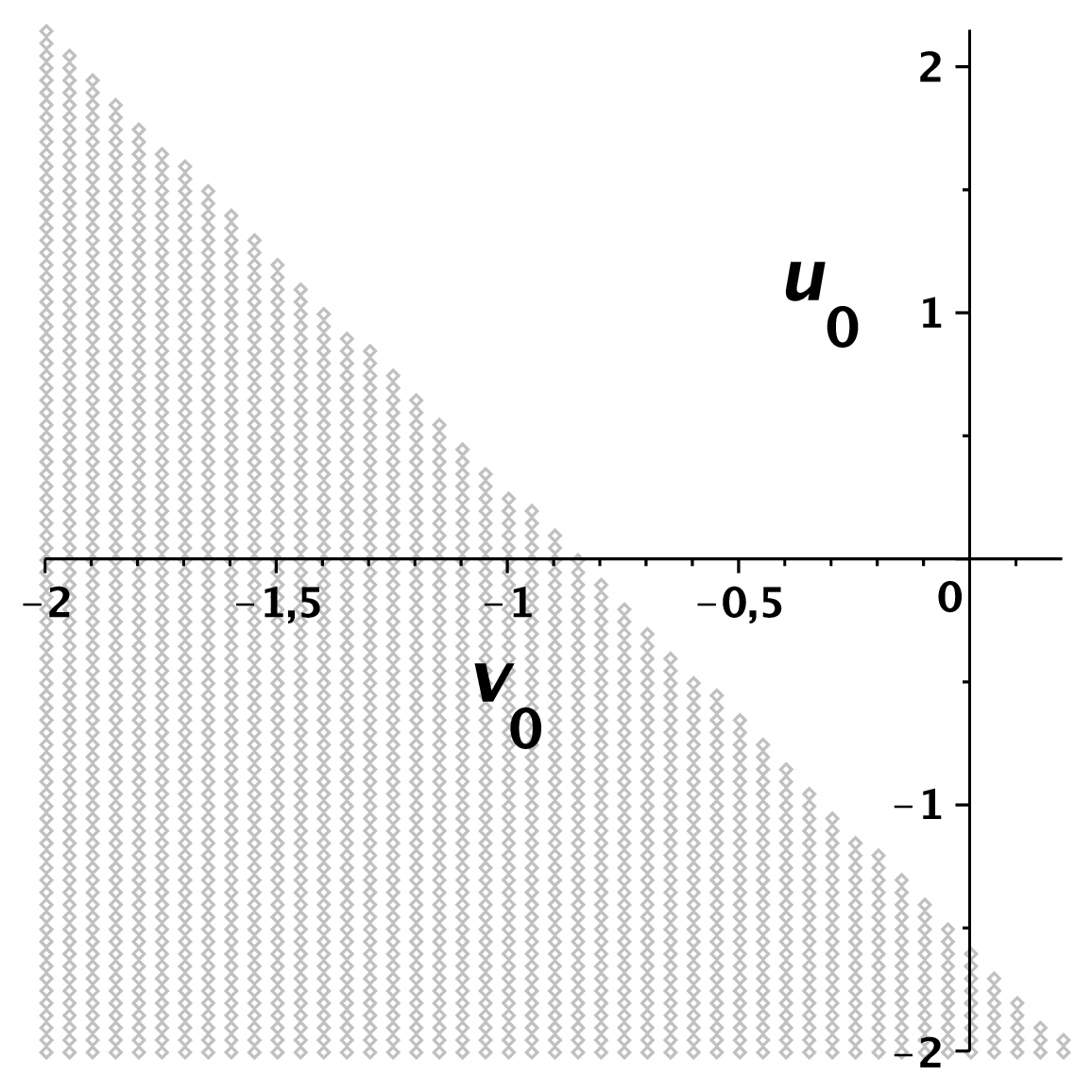}
\end{minipage}
\caption{The structure of the set, corresponding to a singularity formation (dash), for $d=3$, $k=-1$, $c=1$. Left: $F_0=u_0=0$ ($\bv_0=0$). Center: $G_0=v_0=0$ ($\bE_0=0$). Right: $F_0=0$, $G_0=0.2$.}\label{Pic3}
\end{figure}

\section{Discussion}\label{Sec10}
In this paper we show that the method of linearization   can be successfully applied to construct criteria for singularity formation  for some classes of problems related to non-relativistic Euler-Poisson equations without pressure. We show that for some types of initial data and some spatial dimensions, analytical criteria can be obtained, at least in terms of special functions. In all other cases, a criterion for singularity formation is obtained in terms of some auxiliary function. The criterion  can be easily realized numerically. These results can also be useful for estimates that allow one to obtain various types of sufficient conditions for a  singularity formation, for example, in the case of dissipation \cite{RD24}.
 Of course, this is an interesting, but purely mathematical problem, and the benefit from detailed results can be extracted not so much by physicists as by specialists in numerical methods for testing high-precision algorithms \cite{CH18}.
Physicists would be more interested in the relativistic case and the possibility of adding the magnetic field to the problem.
Note that in the context of cold plasma oscillations, the method considered here also provides a criterion for the occurrence of a singularity and gives analytical results for the Davidson model  with the magnetic field \cite{RD_interplay} and for relativistic plasma. However, for the relativistic case, globally smooth solutions to the Cauchy problem do not exist in general (at least in the context of cold plasma), so it is impossible to talk about finding a criterion. However, the question of the lifetime of a smooth solution remains, and the answer can also be obtained in terms of an auxiliary function (based on the proposed numerical algorithm).

Note that the repulsive case with a non-negative background, for which the solutions are oscillating, is rather unpromising from the point of view of non-constant smooth solutions. For example, even in the simplest one-dimensional case, any deviation of the density background from a constant destroys a globally smooth solution \cite{Roz_doping}. The situation is quite different in the attractive case. As was recently shown, in such a situation it is possible to construct globally smooth solutions for non-constant density profiles \cite{Tadmor24}.

In addition, one can consider a wide class of non-strictly hyperbolic equations, generally speaking, not having a physical nature, for which the described technique is applicable. For the case of two equations and spatial dimension one, the problem of finding a criterion for the singularity formation of a solution to the Cauchy problem in terms of the initial data is solved in \cite{T24}.

\section*{Acknowledgements}

 Supported by Russian Science Foundation  grant 23-11-00056 through RUDN University.


\begin{thebibliography}{99}




\bibitem{Bhat23}    M. Bhatnagar, H. Liu,
A complete characterization of sharp thresholds to
spherically symmetric multidimensional pressureless
Euler-Poisson systems, arXiv:2302.04428 (2023).

\bibitem{Bellman} R. Bellman,	Stability theory of differential equations, Dover Books on Mathematics, Courier Corporation (2013).

\bibitem{Brenner}
M. P. Brenner, T. P. Witelski, On spherically symmetric gravitational collapse, J. Stat. Phys., {\bf 93}
863 -- 899 (1998).

\bibitem{Brunelli}
J. C. Brunelli, A. Das, 
On an integrable hierarchy derived from the isentropic gas dynamics, J. Math.
Phys. {\bf 45} (7) 2633 -- 2645 (2004).

\bibitem{Carrillo}J.A. Carrillo, R. Shu,  Existence of radial global smooth solutions to the pressureless Euler-Poisson equations with quadratic confinement. Arch. Rational Mech. Anal. {\bf 247}, 73 (2023).



\bibitem{CH18}  E.V. Chizhonkov, { Mathematical aspects of modelling oscillations and wake waves in plasma}, CRC Press, 2019.



\bibitem{Tadmor24}
Y.-P. Choi, D. Kim, D. Koo, E. Tadmor.
Critical thresholds in pressureless Euler-Poisson
equations with background states, arXiv:2402.12839 [math.AP] (2024).



\bibitem{RD24} M.I.Delova,  O.S. Rozanova, {On radially symmetric oscillations of a collisional cold plasma}, Mathematical Methods in the Applied Sciences, {\bf 47} (11) 8385--8399 (2024).

 \bibitem{RD_interplay} M.I.Delova,  O.S. Rozanova, {The interplay of regularizing factors in the model of upper hybrid oscillations of cold plasma}, Journal of Mathematical Analysis and Applications, {\bf 515} (2)  126449 (2022).





\bibitem{Riccati} G. Freiling,
A survey of nonsymmetric Riccati equations, Linear Algebra and its
Applications 351-352, 243-270 (2002).





\bibitem{ELT}S.Engelberg,
H.Liu, E.Tadmor, Critical thresholds in Euler-Poisson equations,
Indiana University Mathematics Journal, {\bf 50}, 109-157 (2001).

\bibitem{hyper}A. B.
Olde Daalhuis,  Hypergeometric function, in:  NIST Handbook of Mathematical Functions, Cambridge University Press (2010).






\bibitem{Radon} {W. T. Reid,} Riccati differential equations, Academic Press, New York, 1972.

\bibitem{Romanovski} V.G.Romanovski, D.S.Shafer, 
{\em The center and cyclicity problems: A computational
Algebra Approach}, Boston: Birkhauser (2009).

\bibitem{RChZAMP21}
O.S. Rozanova, E.V. Chizhonkov, { On the conditions for the breaking
of oscillations in a cold plasma, { Z. Angew. Math. Phys.},
 }{\bf 72} (2021), 13. 






\bibitem{R22_Rad}
O.S. Rozanova, { On the behavior of multidimensional radially symmetric solutions of the repulsive Euler-Poisson equations}, Physica D: Nonlinear Phenomena  {\bf 443}, 133578 (2023). 

\bibitem{Roz_doping}  O.S. Rozanova, The repulsive Euler-Poisson equations with  variable doping profile, Physica D: Nonlinear Phenomena
{\bf  472} (2) 134454 (2025).


\bibitem{R_exept}  O.S. Rozanova, Criterion of singularity formation for radial solutions of the pressureless Euler-Poisson equations in exceptional dimension, arXiv:2408.13794, submitted.


\bibitem{RT} O.S. Rozanova, M.K.Turzynsky, {On the properties of affine solutions of cold plasma equations}, Communications in Mathematical Sciences, {\bf 22} (1)  215-226 (2024).




\bibitem{Sabatini}
M. Sabatini,  On the period function of Li\'enard systems. J. Differ. Equ. {\bf 152}, 467-487 (1999).

\bibitem{Heun} B. D. Sleeman, V. B. Kuznetzov,   Heun function, in:  NIST Handbook of Mathematical Functions, Cambridge University Press (2010).

\bibitem{Tan}  C. Tan, Eulerian dynamics in multidimensions with radial symmetry. SIAM Journal on Mathematical Analysis,
{\bf 53} (3), 3040-3071 (2021).

\bibitem{T24}  M.K.Turzynsky, {Nonstrictly hyperbolic systems and their application to study of Euler-Poisson equations}, Siberian
Electronic Mathematical Reports, {\bf 21} (2)  215-226 (2024) (arXiv:2410.04597v1).
\end{thebibliography}
\end{document}